\documentclass[preprint, review]{elsarticle}
\usepackage{amsmath}
\usepackage{amssymb}
\usepackage{amsthm}
\usepackage{mathrsfs}
\usepackage[nodots]{numcompress}

\journal{Nonlinear Analysis: Real World Applications}

\theoremstyle{plain} %default

\newtheorem{lem}{Lemma}

\newtheorem{thm}{Theorem}

\theoremstyle{definition}

\theoremstyle{remark}

\newcommand{\pd}{\,\partial}
\newcommand{\comment}[1]{}

\begin{document}

\begin{frontmatter}



\title{Equivalence transformations of  Euler-Bernoulli equation}

\author{J.C. Ndogmo\corref{cor1}}
\ead{jean-claude.ndogmo@wits.ac.za}

\cortext[cor1]{Tel: $+$27 11 717 6237; $\quad$ Fax: +27 86 579 6248}

\address{School of Mathematics,
University of the
Witwatersrand,
Private Bag 3, Wits 2050,
South Africa}

\begin{abstract}
 We give  a
determination of the equivalence group of Euler-Bernoulli equation
and  of one of its generalizations, and thus derive some symmetry properties
of this equation.
\end{abstract}

\begin{keyword}
equivalence group \sep  arbitrary functions \sep
Euler-Bernoulli equation \sep  differential operators \sep symmetry group

\MSC[2008] 58J70 \sep  34C20 \sep  35A30
\end{keyword}

\end{frontmatter}

\section{Introduction}
\label{s:intro}

Consider a collection $\mathcal{F}$ of differential equations of the form
\begin{equation}
\label{eq:geq} \Delta (x, y_{(n)}; \mathscr{A} )=0,
\end{equation}
where $x= (x^1, \dots, x^p)$ is the set of independent variables
and $y_{(n)}$ denotes $y$ and all its derivatives up to the order
$n,$ while the parameter $A$ denotes collectively the set of all arbitrary
functions specifying the family element in
$\mathcal{F}.$ On the other hand, let $G$ be the Lie pseudo-group  of point transformations  of the
form
\begin{equation}\label{eq:gtransfo}
x= \varphi (z, w), \qquad y = \psi (z, w),
\end{equation}
where  $z= (z^1, \dots, z^p)$ is the new set of independent
variables, and $w=w(z)$ is the new dependent variable. As Tresse explained in
\cite[P. 11]{tresse}, the elements of $G$ depend in general on arbitrary functions and not on arbitrary constants, and $G$ is  therefore infinite-dimensional. We say that $G$ is the group of equivalence transformations of \eqref{eq:geq} if it is the largest Lie
pseudo-group of point transformations that map \eqref{eq:geq} into
an equation of the same form, that is, if in terms of the same
function $\Delta$ appearing in \eqref{eq:geq}, the transformed equation has an expression of the form
\begin{equation} \label{eq:trfoeq}
\Delta (z, w_{(n)}; B )=0,
\end{equation}
where $B$ denotes collectively the new set of arbitrary functions of the equation. Equivalence groups play an important role in, amongst others,  the classification of differential equations by means of their invariant functions, and  they are also a valuable tool for the identification of these invariant functions \cite{schw-pap, schw-bk, olv-fls2, olv-maurer2}. \par

 In this paper we obtain the equivalence transformations of Euler-Bernoulli equation (\textsc{EB}). The difficulty with this equation is due to the high order of derivatives that it involves,  especially in the two-dimensional case, and this  makes the standard method of finding the equivalence transformations by direct analysis intractable. We therefore recourse to the transformation of differential operators that we discuss in the next section. Next, we proceed to the determination of the equivalence group of (\textsc{EB}) and of one of its generalizations, and derive some symmetry properties of this equation. This includes the full symmetry group of the equation, which was obtained only partially in \cite{waf-eub1} in terms of some undetermined set of functions.

\section{Transformation of differential operators}
  An expression of differential operators in terms of new variables
yields an effective algorithm for implementing the transformation of
a differential equation under a given change of variables, and this
can also greatly simplify calculations, based on the properties of
these differential operators.\par

Suppose for instance that the general change of variables
\eqref{eq:gtransfo} that maps \eqref{eq:geq} into \eqref{eq:trfoeq}
is given more explicitly in the form
\begin{subequations}\label{eq:gnchg2}
\begin{align}
x^i &= \varphi_i (z, w), \quad z= (z^1, \dots, z^p),\quad w= w(z) \\
u   &= \psi(z, w) \equiv T(z).
\end{align}
\end{subequations}
 Then the last equality gives
\begin{equation}\label{eq:linsyst1}
\frac{\pd\, T(z)}{ \pd z^i} = \sum_j u_{x^j} \frac{\pd
\varphi_j}{\pd z^i},\quad \text { for $i=1, \dots, p$,}
\end{equation}
and for invertible transformations \eqref{eq:gnchg2}, the linear
system \eqref{eq:linsyst1} yields solutions of the form
$$
u_{x^j} = \sum_{i=1}^p \Psi_{\!j}^i \frac{\pd\, T(z)}{ \pd z^i}
,\quad \text { for $j=1, \dots, p$},
$$
for some functions $\Psi_{\!j}^i= \Psi_{\!j}^i(z).$ In other words,
an expression of differential operators in terms of the new
independent variables is given by
\begin{equation}\label{eq:linop1}
\pd_{x^j} = \Psi_{\!j}, \quad \text{ where }\quad \Psi_{\!j}= \sum_i
\Psi_{\!j}^i  \pd_{z^i},
\end{equation}
We shall make use of some properties of the linear differential
operators $\Psi_{\!j}$ to simplify calculations, especially when
dealing with high-order derivatives.

\section{Determination of the equivalence group}
  We now move on to consider in this section the problem of finding
the equivalence group  for the one-dimensional Euler-Bernoulli
equation, which is a model for calculating the load-carrying and
deflection characteristics of beams. Since the late $19$th century,
following the successful demonstration of this theory in the
constructions of the Eiffel Tower and the Ferris wheel,
Euler-Bernoulli equation (also known as Engineer's beam theory)
became the cornerstone of engineering.\par

This equation can be put in the form
\begin{equation}\label{eq:eub}
\frac{\pd^2}{\pd x^2} \left( f(x) \frac{\pd^2 u}{\pd x^2}\right) +
m(x) \frac{\pd^2 u}{\pd t^2} =0, \qquad t>0, \quad 0<x< \mathcal{D},
\end{equation}
where $f(x)>0$ is the flexural rigidity, $m(x)>0$ is the lineal mass
density and $u(t,x)$ is the transversal displacement  at time $t$
and position $x$ from one end of the beam taken as origin.

In the relatively recent research literature,
Euler-Bernoulli equation, which is often simply referred to as
Euler-Bernoulli beam or beam equation,  has been frequently
discussed. In particular, Gotllieb \cite{gottlieb} has investigated
the iso-spectral properties of this equation and its non-homogeneous
variants in connection with the unit beam (i.e. $f=1$ and $m=1$).
More recently, Wafo Soh \cite{waf-eub1} considered the {\em
equivalence problem} for Euler-Bernoulli beam from the Lie symmetry
approach, while Morozov and Wafo Soh \cite{waf-mrz1} investigated
the same problem using Cartan's equivalence method. \par

It is clear that in any change of variables that preserves the form
of \eqref{eq:eub}, the transformation of $x$ should depend on a
single variable (due to the form of $f$ and $m$), while the
linearity of the equation forces the transformation of $u$ to be
linear, and that for the independent variables not to involve the
dependent variable. We are therefore led to look for equivalence
transformations of \eqref{eq:eub} in the form
\begin{equation}\label{eq:eqveub1}
t= R(y,z), \quad x= S(z), \quad u = L(y,z)w+ J(y,z)
\end{equation}

Under \eqref{eq:eqveub1}, the coefficient $\gamma_1$ of $w_{yyyy}$
in the transformed equation takes the form
$$
\gamma_1 = (f R_z^4 L)/ (R_y^4 S_z^4)
$$
and its vanishing  implies that $R= R(y).$
Inserting this new expression for $R$ into \eqref{eq:eqveub1}, the
vanishing of the coefficient of $w_y$ in the new transformed
equation yields the condition
$$
2 R_y L_y - L R_{yy} =0,
$$
and thus $L= h(z) \sqrt{R_y},$ for some function $h.$ If we now
transform \eqref{eq:eub} using this additional information and write
down the  coefficient $\gamma_2$ of $w_z$ as a polynomial in $f(S)$
and its derivatives, then the coefficient $\delta_1$ of $f_{S,S}$ in
this expansion is
$$
\delta_1 = S_z^4 (2 h_z S_z - h S_{zz}),
$$
and its vanishing gives $h= k_1 \sqrt{S_z},$ where $k_1$ is a
constant, so that $u= k_1 (R_y S_z)^{1/2} w + J$. With this
expression for $h,$ the vanishing of the coefficient of $f_S$ in
$\gamma_2$ leads to the condition
$$
-3 \,S_{zz}^2 + 2 \,S_z S_{zzz}=0,
$$
and this in turn shows that
$$
x \equiv S= (k_2 z+ k_3)/(k_4 z+ 1), \qquad k_2- k_3 k_4 \neq 0,
$$
for some constants $k_2, \dots, k_4.$ With all the information
obtained up to this point on $R, S,$ and $L,$ the coefficient
$\gamma_3$ of $w$ now takes the form
$$
\gamma_3 =  -\frac{4 (k_2-k_3 k_4)^8 m \left(3 R_{yy}^2-2 R_y
R_{yyy}\right)}{(1+k_4 z)^{16}},
$$
showing that $R$ must satisfy an equation similar to that for $S,$
namely
$$
3 R_{yy}^2 - 2 R_y R_{yyy}=0.
$$
Consequently, we must  have
\begin{equation}\label{eq:R1}
t \equiv R= (k_{\,5} y + k_6)/(k_7 y+1), \qquad k_{\,5} - k_6 k_7
\neq 0,
\end{equation}
for some arbitrary constants $k_{\,5}, \dots, k_7.$ With this new
information and all others obtained thus far, \eqref{eq:eqveub1} now
takes the form
\begin{subequations} \label{eq:eqveub2}
\begin{align}
t &= \frac{k_{\,5} y+ k_6}{k_7 z+1}, \quad x= \frac{k_2 z+ k_3}{k_4
z+1}  \label{eq:eqveub2a}\\
 u &= k_1\frac{(k_2 - k_3 k_4)^{1/2}(k_5- k_6
k_7)^{1/2}}{(1+ k_7 y)(1+ k_4 z)}\, w + J(y,z).
\end{align}
\end{subequations}

Under \eqref{eq:eqveub2}, the transformation of  Euler-Bernoulli
beam retains a constant component, that is, a component (or sum of
terms) free from $w$ and its derivatives, and which appears to
contain the only terms involving $J$ and its derivatives. In fact
this component is simply the transformation of $u=J$ under the
change of independent variables \eqref{eq:eqveub2a}, and it is
readily found that the vanishing of this constant component amounts
to setting
\begin{equation}\label{eq:J1}
 J(y,z)= \frac{k_8-k_{10} k_7^2 y+k_4 \left(-k_9+k_{11} k_7^2
y\right)+k_4^2 \left(-k_9 z+k_{11} k_7^2 y z\right)}{k_4 k_7 (1+k_7
y) (1+k_4 z)}
\end{equation}
where $k_8, \dots k_{11}$ are some additional constants of
integration.

If we now perform a transformation of \eqref{eq:eub} under the
change of variables \eqref{eq:eqveub2} with $J$ as in \eqref{eq:J1},
the transformed equation takes on the form

\begin{subequations}\label{eq:eubytr}
\begin{align}
\quad& \frac{\pd^2}{\pd z^2} \left( F(y,z) \frac{\pd^2 w}{\pd
z^2}\right) + M(y,z) \frac{\pd^2 w}{\pd y^2} =0,\\
\intertext{ where }
F(y,z) &= \frac{(1+k_4 z)^5}{(1+k_7 y) (1+k_4 z)} f(S), \\
M(y,z) &= \frac{(1+k_7 y)^3 (1+k_4 z)^3(k_2-k_3 k_4)^4 }{(k_5-k_6
k_7)^2 (1+k_4 z)^7}\, m(S).
\end{align}
\end{subequations}

Therefore, the change of variables given by \eqref{eq:eqveub2} and
\eqref{eq:J1} will represent an equivalence transformation of
\eqref{eq:eub} if and only if the functions $F$ and $M$ in
\eqref{eq:eubytr} do not depend on $y,$ and this is possible if and
only if $k_7=0.$ However, as in \eqref{eq:J1} the value $k_7=0$ is
not allowed, we must recompute the appropriate expression for $J$
when $k_7=0.$ So, letting $k_7=0$ in \eqref{eq:R1} and finding the
 corresponding value for $J$ gives
\begin{equation}\label{eq:J2}
J=\frac{-k_0-k_9 y+k_4 (k_8+k_{10} y)+k_4^2 (k_8+k_{10} y) z}{k_4
(1+k_4 z)},
\end{equation}
for some new constants $k_0,$ and $k_8, \dots, k_{10},$ and the
transformed equation now takes on the required form
\begin{subequations}\label{eq:eubnor}
\begin{align}
\quad& \frac{\pd^2}{\pd z^2} \left( F(z) \frac{\pd^2 w}{\pd
z^2}\right) + M(z) \frac{\pd^2 w}{\pd y^2} =0,\\
\intertext{ where }
F(z) &=  (1+k_4 z)^5 \left[\frac{ (k_2- k_3 k_4) k_5}{(1+k_4 z)^2}\right]^{1/2}f(S), \\
M(z) &= \frac{(k_2- k_3 k_4)^5 }{k_5^2 (1+k_4 z)^3}
\left[\frac{(k_2-k_3 k_4) k_5}{(1+k_4 z)^2} \right]^{1/2} m(S).
\end{align}
\end{subequations}
We have thus obtained the following result.
\begin{thm}\label{th:eqveub}
The group  $G$ of equivalence transformations of Euler-Bernoulli
equation \eqref{eq:eub} is given by the equivalence transformations
\begin{subequations} \label{eq:eqveubfn}
\begin{align}
t &=k_{\,5} y+ k_6, \quad x= \frac{k_2 z+ k_3}{k_4
z+1}\\
 u &= k_1 \left(\frac{(k_2 - k_3 k_4)k_5}{(1+ k_4 z)^2}\right)^{1/2}   \, w +
 J,\\
 \intertext{where}
 J &= \frac{-k_0-k_9 y+k_4 (k_8+k_{10} y)+k_4^2 (k_8+k_{10} y) z}{k_4
(1+k_4 z)}, \notag
\end{align}
\end{subequations}
is the function obtained in \eqref{eq:J2}, and the transformed
equation is given by \eqref{eq:eubnor}.
\end{thm}

  We also note that in the successive transformations of \eqref{eq:eub}
which led to \eqref{eq:eubytr}, the arguments of the arbitrary
function $f$ and $m$ were treated as dummy, by being ignored. Thus,
in view of \eqref{eq:eubytr}, we readily obtain the following result
for an extended  version of Euler-Bernoulli Equation.
\begin{lem}\label{le:eqveubt}
For the generalized version of Euler-Bernoulli equation of the form
\begin{equation}\label{eq:eubt}
\frac{\pd^2}{\pd x^2} \left( f(t, x) \frac{\pd^2 u}{\pd x^2}\right)
+ m(t,x) \frac{\pd^2 u}{\pd t^2} =0,
\end{equation}
in which  $f$ and $m$ are assumed to be functions of both $t$ and
$x,$ an {\em equivalence subgroup} is given by the equivalence
transformations  \eqref{eq:eqveub2} and \eqref{eq:J1}, and the
transformed equation takes on the form \eqref{eq:eubytr}.
\end{lem}
 We may not however assert as yet that \eqref{eq:eqveub2}
and \eqref{eq:J1} define the equivalence group of \eqref{eq:eubt},
because due to the form of $f$ and $m$ in \eqref{eq:eub}, these
transformations were obtained under the assumption that $x=S(z),$
depends only on a single variable. But this restriction in no longer
{\em a priori} permitted in the case of \eqref{eq:eubt} where $f$
and $m$ are functions of both $t$ and $x.$ \par
The equivalence group of \eqref{eq:eubt} should thus be looked for
in the form
\begin{equation}\label{eq:treubt1}
t= R(y,z), \quad x=S(y,z),\quad \text{ and } u= L(y,z) w+ J(x,y),
\end{equation}
where as usual, the functions $R, S, L$ and $J$ are to be found. It
first follows from the lemma  and the invertibility of
\eqref{eq:treubt1}, that
\begin{subequations}\label{eq:jacoeubt}
\begin{align}
R_y S_z L & \neq 0, \qquad  \varpi L \neq 0, \vspace{-5.5mm}\\
\intertext{$\vspace{-3.5mm}$ where }
\varpi &= R_y S_z - R_z S_y.
\end{align}
\end{subequations}

 It is also clear that in the transformed
equation, terms of the highest order may only come from the
transformation of $u_{xxxx}.$ Moreover, to find such terms, we may
assume that $u= L\, w$ in \eqref{eq:treubt1}, and since $L \neq 0,$
by Leibnitz rule we may assume for such purpose that $u=w(y,z).$
Since the term $( R_z^4 /\varpi)w_{yyyy}$ appears in the polynomial
expansion of the transformation of $u_{xxxx},$ the vanishing of this
term yields  $R=R(y).$\par
On the other hand, it follows from the formulas \eqref{eq:linop1}
that
\begin{equation}\label{eq:opeubt11}
\pd_t = \frac{1}{\varpi} (S_z \pd_y - S_y \pd_z), \quad \pd_x=
\frac{1}{\varpi} (-R_z \pd_y + R_y \pd_z),
\end{equation}
 and for $R=R(y),$ this reduces to
\begin{equation}\label{eq:opeubt11}
\pd_t = \frac{1}{R_y}  \pd_y - \frac{S_y}{R_y S_z} \pd_z, \qquad
 \qquad \pd_x= \frac{1}{S_z} \pd_z.
\end{equation}
Consequently,  terms involving derivatives of $w$ w.r.t. $y$ in the
transformed equation can only come from the transformations of
derivatives of $u$ w.r.t. $t,$ i.e. from the transformation of
$u_{tt}.$ So, transforming $u_{tt}$ by setting $u= L(y,z) w$ in
\eqref{eq:treubt1} shows that the term in $w_{yz}$ has coefficient

$$
-(2 S_y L)/ (R_y^2 S_z).
$$
The vanishing of this term gives  $S=S(z),$ and hence the
transformations of \eqref{eq:eubt} must also be sought in the form
\begin{equation}\label{eq:treubt22}
t= R(y), \quad x= S(z), \quad u = L(y,z)u+ J(y,z),
\end{equation}
just as in the case of  \eqref{eq:eub}. Given that in the
transformations of \eqref{eq:eub} starting with a change of
variables of the form \eqref{eq:treubt22}, and leading to
\eqref{eq:eqveub2}, \eqref{eq:J1}, and \eqref{eq:eubytr}, the
arguments of the arbitrary functions $f$ and $m$ are ignored, and
because \eqref{eq:eub} and \eqref{eq:eubt} differ only by these
arguments, we may  conclude that the equivalence transformations of
\eqref{eq:eubt} must also be of the form \eqref{eq:eqveub2}-
\eqref{eq:J1}. Since the resulting transformations of
\eqref{eq:eubt} which are also given by \eqref{eq:eubytr} are of the
required form, we have thus obtained the following result.
\begin{thm}
The equivalence transformations
\begin{align*}
t &= \frac{k_{\,5} y+ k_6}{k_7 z+1}, \quad x= \frac{k_2 z+ k_3}{k_4
z+1}  \label{eq:eqveub2a}\\
 u &= k_1\frac{(k_2 - k_3 k_4)^{1/2}(k_5- k_6
k_7)^{1/2}}{(1+ k_7 y)(1+ k_4 z)}\, w + J(y,z)\\
\intertext{where}
 J(y,z)&= \frac{k_8-k_{10} k_7^2 y+k_4 \left(-k_9+k_{11} k_7^2
y\right)+k_4^2 \left(-k_9 z+k_{11} k_7^2 y z\right)}{k_4 k_7 (1+k_7
y) (1+k_4 z)},
\end{align*}

and which are given by \eqref{eq:eqveub2} and \eqref{eq:J1}, define
the equivalence group of the extended form \eqref{eq:eubt} of
\eqref{eq:eub}.
\end{thm}

%%%%%%%%%%%%%%%%%%%%%%
%%%%%%%%%%%%%%%%%%%%%%%%%%%%%%
%%%%%%%%%%%%%%%%%%%%%%
\section{Symmetry properties}
\label{s:sym}
 It is clear that the equivalence group of a differential equation
contains according to our definition the largest symmetry group
$G_e$ of the equation that is free from the arbitrary functions of
the equation. However, such a symmetry group is only a subgroup of
the so-called principal Lie group, which is the largest symmetry
group that holds for all arbitrary functions, but which may depend
on these arbitrary functions.\par

Symmetries of Euler-Bernoulli equation were
calculated in \cite{waf-eub1}, but were obtained only in a
conditional form, and under the assumption that the functions $f$
and $m$ satisfy certain unsolved complicated equations.\par

It appears from the form of \eqref{eq:eubnor} that
\eqref{eq:eqveubfn} will represent a symmetry transformation of
\eqref{eq:eub} only if the transformation of $x$ is the identity
transformation, $x=z.$ This means that we have to assume $k_4=0$ in
\eqref{eq:eqveubfn}, and because the corresponding expression for
the denominator of $J$ vanishes for  $k_4=0,$  we have to recompute
the appropriate value of $J.$ To do so, we transform \eqref{eq:eub}
using \eqref{eq:eqveub2}, with $k_7=k_4= k_3=0,$ and $k_2=1.$ It
then follows that \eqref{eq:eub} will be invariant under these
transformations if and only if $k_5=1,$ and
\begin{equation}\label{eq:J3}
J= y p_2+ p_4 + z \left(y p_1+p_3\right).
\end{equation}
We thus have the following result.
\begin{thm}\label{th:Ge}
 The largest symmetry group of Euler-Bernoulli equation that
 is independent of its arbitrary functions is a six-parameter group $G_e$ given by the
 transformations
 \begin{equation}\label{eq:gpGe}
t= y+ p_{6},\quad x= z,\quad \text{ and } \quad u= p_{\,5} w + J,
 \end{equation}
where $J$ is given by \eqref{eq:J3}, and with corresponding symmetry
generator
 \begin{equation}\label{eq:symGe}
\mathbf{v} = p_6 \pd_t + (p_4 + p_2 t+ + p_3 x + p_1 t x+ p_{\,5}
u)\pd_u,
 \end{equation}
where $p_1, \dots, p_6$ are the group parameters.
\end{thm}
We note that the function $J$ in \eqref{eq:J3} is the fundamental
solution of  Euler-Bernoulli equation that does not depend on its
arbitrary functions. Owing to the linearity of \eqref{eq:eub}, if we
replace
 the function $J$ in \eqref{eq:gpGe} by any other solution of
 \eqref{eq:eub},  the resulting transformation remains a symmetry \cite{olv-bk2, bluman},
 but which may however depend on the arbitrary functions of
the equation, and this is contrary to our assumptions, and for
practical considerations such solutions are not always available.\par

On the other hand, due to the superposition   principle of linear equations,  \eqref{eq:gpGe}  with $J$
replaced by an arbitrary solution $S$ of \eqref{eq:eub}, remains a symmetry of Euler-Bernoulli equation, and one readily verifies that this is the most general symmetry group of this equation.

%\section{Concluding remarks}

%%%%%%%%%%%%%%%%%%%%%%%%%%%%%%%%%%%%%%%%%%%%%%%%%%%%%%%%%%%%%%%%%%%%%%%%%%%
%%%%%%%%%%%%%%%%%%%%%%%%%%%%
%%%%%%%%%%%%%%%%%%%%%%%%%%%%%%%%%%%%%%%%%%%%%%%%%%%%%%%%%%%%%%%%%%%%%%%%%%%

\bibliographystyle{model1-num-names}

\end{document}